\documentstyle[aps]{revtex}   

\begin{document}
\draft
\title{ The Local Moduli of Sasakian 3-Manifolds }
\author{Brendan S. Guilfoyle \footnote{Email: brendan.guilfoyle@ittralee.ie} }
\address{Mathematics Department, Institute of Technology Tralee, Tralee, Co. Kerry, Ireland.}
\date{\today}
\maketitle
\begin{abstract}
The Newman-Penrose-Perjes formalism is applied to Sasakian 3-manifolds and the local form of the metric and contact structure is presented.  The local moduli space can be parameterised by  a single function  of two variables and it is shown that, given any smooth function of two variables, there exists locally  a Sasakian structure with scalar curvature equal to this function.  The case where the scalar curvature is constant ($\eta$-Einstein Sasakian metrics) is completely solved locally.  The resulting Sasakian manifolds include $S^3$, $Nil$ and $\tilde{SL_2R}$, as well as the Berger spheres.  It is also shown that a conformally flat Sasakian 3-manifold is Einstein of positive scalar curvature.
\end{abstract}

\newcommand{\nc}{\newcommand}
\nc{\sigf}{\Sigma_\phi}
\nc{\bg}{\begin{eqnarray}}
\nc{\ed}{\end{eqnarray}}
\nc{\bp}{\bar{\Psi}}
\nc{\bm}{\bar{m}}
\nc{\ul}{\underline}
\nc{\bee}{\begin{equation}}
\nc{\ee}{\end{equation}}
\nc{\beer}{\begin{eqnarray}}
\nc{\eer}{\end{eqnarray}}
\nc{\dv}{\textstyle\frac}
\nc{\sdv}{\scriptstyle\frac}
\nc{\hf}{\dv{1}{2}}
\nc{\qtr}{\dv{1}{4}}
\nc{\shf}{\sdv{1}{2}}
\nc{\nn}{\nonumber}

\section{Introduction}
Sasakian geometry  consists of a riemannian metric and a contact structure adapted to the metric in a certain sense.  First introduced in 1960 \cite{sas}, Sasakian geometry is the odd-dimensional equivalent of K\"{a}hler geometry.   For review papers see Boyer and Galicki \cite{bag1} \cite{bag2}.

Geiges \cite{geiges} has classified the topology of all closed Sasakian 3-manifolds:

\vspace{0.2in}

\noindent {\bf Theorem}

{\it A closed 3-manifold admits a Sasakian structure iff it is diffeomorphic to one of the following:
\begin{description}
\item[(a)] $S^3$
\item[(b)] $\tilde{SL_2R}$
\item[(c)] $Nil$
\end{description}
or a quotient by a discrete subgroup of the appropriate  isometry group acting freely.
}

\vspace{0.2in}

The proof of this theorem rests on the fact that if $M$ is a closed Sasakian 3-manifold, then $M\times R$  is a locally conformally K\"{a}hler manifold with parallel Lee form.  
In addition, Belgun \cite{belgun}  has classified the CR structures associated with closed non-spherical Sasakian 3-manifolds.  It has also been shown that any Sasakian-Einstein structure on $R^3$ is diffeomorphic to the standard Sasakian structure on an open subset of $S^3$ \cite{guil}. 

The purpose of this paper is to complete the local classification of Sasakian structures.  Our main theorem is

\vspace{0.2in}

\noindent {\bf Theorem}
{\it The local space of smooth Sasakian structures is given by a real function of two variables.
}

\vspace{0.2in}

In essence, the structure theorem says that locally, all Sasakian structures arise as a Kaluza-Klein construction on a line bundle over an open set $U\subset R^2$, the real function mentioned above being the conformal factor of the riemann metric on $U$, when written in isothermal coordinates.  The Sasakian structures are in one to one correspondence with the riemann structures on $U$.  In addition, all Sasakian $\eta$-Einstein structures are found.  These include $S^3$, $Nil$ and $\tilde{SL_2R}$, as well as the Berger spheres. Many of these are examples of Sasakian 3-manifolds with negative scalar curvature.  

The main tool used here is the Newman-Penrose-Perjes spin coefficient formalism, details of which can be found in \cite{guil}.  This technique is also used to show that there are no non-Einstein conformally flat Sasakian 3-manifolds, a theorem first proved by Tanno \cite{tanno2}.  

This paper is organized as follows.  In the next Section the Newman-Penrose-Perjes formalism for riemannian 3-manifolds is outlined.  In Section 3 Sasakian structures are introduced and  are characterized in terms of the divergence, shear and twist of the associated geodesic foliation.  The main theorem for the moduli space of Sasakian structures is also proven.  Section 4 addresses Sasakian $\eta$-Einstein metrics and finds that there are three families of these, each with examples of (constant) negative curvature.  In Section 5, we give a short proof of the fact that the only conformally flat Sasakian 3-manifolds are of constant curvature.

\section{The Newman-Penrose-Perjes Formalism }

Let ($M$, $g$) be a riemannian 3-manifold and consider an orthonormal frame $\{e_{0}, e_1, e_2\}$, and dual basis of 1-forms $\{\theta^{0}, \theta^1, \theta^2\}$.  Introduce a complex frame $\{e_{0}, e_+, e_-\}$ by 
\[
e_+=\frac{1}{\sqrt{2}}\left(e_1-i\;e_2\right) \qquad\qquad e_-=\frac{1}{\sqrt{2}}\left(e_1+i\;e_2\right),
\]
with dual basis of 1-forms $\{\theta^{0}, \theta^+, \theta^-\}$.

We define the complex {\it spin coefficients} by
\[
\gamma_{mnp}=\nabla_je_{mi}\;e_n^ie_p^j,
\]
where $\nabla$ is the covariant derivative associated with $g$ and the indices $m$, $n$, $p$ range over $0$, $+$, $-$.  Thus
\[
\gamma_{mnp}=-\gamma_{nmp}.
\]
We break covariance and introduce the complex optical scalars
\[
\gamma_{+0-}=\rho \qquad\qquad \gamma_{+0+}=\sigma \qquad\qquad \gamma_{+--}=\tau
\]
\[
\gamma_{+00}=\kappa \qquad\qquad \gamma_{+-0}=\epsilon.
\]
Geometrically, $\kappa=0$ if and only if $e^i_{0}$ is tangent to a geodesic congruence in $M$.  The real part of $\rho$ measures the divergence, the imaginary part of $\rho$ measures the twist and $|\sigma|$ measures the shear of the congruence of geodesics.

We introduce the differential operators
\[
D=e_{0}^i\;\frac{\partial}{\partial x^i} \qquad \qquad \delta=e_+^i\;\frac{\partial}{\partial x^i}  \qquad \qquad \overline{\delta}=e_-^i\;\frac{\partial}{\partial x^i},
\]
and if $f$ is any function, the commutators of these operators work out to be
\[
(D\delta-\delta D)f=[(\overline{\rho}+\epsilon)\delta +\sigma\overline{\delta}+\kappa D]f
\]
\[
(\delta\overline{\delta}-\overline{\delta}\delta)f=[\overline{\tau}\overline{\delta} -\tau\delta+(\overline{\rho}-\rho) D]f.
\]
The components of the Ricci tensor in terms of the spin coefficients are:
\[
R_{00}=D\rho+D\overline{\rho}-\overline{\delta}\kappa-\delta\overline{\kappa}+\tau\kappa+\overline{\tau}\;\overline{\kappa}-2\kappa\overline{\kappa}                       
        -2\sigma\overline{\sigma}-\rho^2-\overline{\rho}^2
\]
\[
R_{++}=-\delta\kappa+D\sigma-2\epsilon\sigma-\overline{\tau}\kappa-\kappa^2
        -\sigma\overline{\rho}-\rho\sigma 
\]
\[
R_{0+}=-\overline{\delta}\sigma+\delta\rho+2\tau\sigma+\kappa\rho-\kappa\overline{\rho} 
\]
\[
R_{0-}=-\overline{\delta}\epsilon+D\tau +\kappa\overline{\sigma}-\rho\overline{\kappa}+\epsilon\tau
        -\epsilon\overline{\kappa}+\overline{\tau}\;\overline{\sigma}-\tau\rho 
\]
\[
R_{+-}= -\overline{\delta}\kappa+D\rho+\delta\tau+\overline{\delta}\;\overline{\tau}+\epsilon\rho
        -\epsilon\overline{\rho}-\kappa\overline{\kappa}+\kappa\tau-\rho\overline{\rho}-\rho^2
        -2\tau\overline{\tau},
\]
while the scalar curvature is
\[
{\shf}R=-2\delta\overline{\kappa}+2D\overline{\rho}+\delta\tau+\overline{\delta}\overline{\tau}-2\kappa\overline{\kappa}+2\overline{\kappa}\;\overline{\tau}-2\overline{\rho}^2-\sigma\overline{\sigma}+\epsilon\rho-\epsilon\overline{\rho}-\rho\overline{\rho}-2\tau\overline{\tau} .
\]
In addition, we have the following identities from the symmetries of the curvature tensor
\bee \label{e:id1}
D\rho-\overline{\delta}\kappa+\kappa\tau-\rho^2 = D\overline{\rho}-\delta\overline{\kappa}
        +\overline{\kappa}\;\overline{\tau}-\overline{\rho}^2
\ee
\bee \label{e:id2}
\delta\overline{\sigma}-\overline{\delta}\overline{\rho}-\overline{\tau}\;\overline{\sigma}-\overline{\kappa}\;\overline{\rho} = 
   \overline{\delta}\epsilon-D\tau -\kappa\overline{\sigma} -\epsilon\tau+\epsilon\overline{\kappa} +\tau\rho.
\ee
and the Bianchi identities
\bee\label{e:bian1}
DE_{00}+\overline{\delta}E_{0+}+\delta E_{0-}+(\rho+\overline{\rho})(E_{+-}-E_{00}) +(2\overline{\kappa}-\tau)E_{0+} +(2\kappa-\overline{\tau})E_{0-} +\overline{\sigma}E_{++}+\sigma E_{--}=0
\ee
\bee\label{e:bian2}
DE_{0+}+\overline{\delta}E_{++}+\delta E_{+-}+\kappa(E_{+-}-E_{00})-(\epsilon+2\rho+\overline{\rho})E_{0+}-\sigma E_{0-}+(\overline{\kappa}-2\tau)E_{++}=0
\ee
where the energy momentum tensor is defined by $E_{ij}=R_{ij}-{\shf}Rg_{ij}$.

\section{Sasakian Structures}

In 1960 Sasaki introduced a type of metric-contact structure which can be thought of as the odd-dimensional version of K\"{a}hler geometry.  There are many equivalent ways of defining such a structure, and we will adopt the following:

A riemannian manifold ($M$, $g$) is {\it Sasakian} if there exists a unit length Killing vector field $e_{0}$ on $M$ so that the riemannian curvature satisfies the condition
\bee\label{e:sas}
R(X,e_{0})Y=g(e_{0},Y)X-g(X,Y)e_{0},
\ee
for any vector fields $X$ and $Y$ on $M$.

We connect this with the spin coefficients by adapting the frame to the Killing vector in the following:

\vspace{0.2in}

\noindent {\bf Proposition 1}

{\it
A unit length vector $e_{0}$ is a Killing vector on a riemannian 3-manifold if and only if it is tangent to a divergencefree, shearfree congruence.  In addition, such a congruence is geodesic.
}

\noindent{\bf Proof:}

Suppose that $e_{0}$ is a unit Killing vector.  That is
\[
e_{0\;i}\;e_{0}^{\;i}=1 \qquad\qquad \nabla_{(i}e_{0\;j)}=0.
\]
The contraction of the second of these says that $e_{0}$ is tangent to a divergencefree congruence.  From the definition of the shear, it is clear then that the congruence is shearfree.

Conversely, suppose that the congruence is divergencefree and shearfree.  Then, from the definition of the shear
\[
0= \nabla_{(i}e_{0\;j)}\;\nabla^{(i}e_0^{\;j)}.
\]
Since $g$ is positive definite, we can conclude that $e_{0}$ is a Killing vector field.  Finally
\[
0=e_{0}^{\;j} \;\nabla_{(j}\;e_{0i)}= e_{0}^{\;j} \;\nabla_{j}\;e_{0i}+e_{0}^{\;j} \;\nabla_{i}\;e_{0j}= e_{0}^{\;j}\;\nabla_{j}\;e_{0i},
\]
since  $e_{0}$ has constant length.  Thus the congruence is geodesic. $\Box$

\vspace{0.2in}

\noindent {\bf Theorem 2 }

{\it
A riemannian 3-manifold ($M$, $g$) is Sasakian if and only if there exists a geodesic congruence which
\begin{description}
\item[(i)] is divergence-free
\item[(ii)] is shear-free
\item[(iii)] has twist equal to one.
\end{description}
}

\noindent{\bf Proof:}

In the spin coefficient notation adapted to the Killing vector field a Sasakian 3-manifold has $\kappa=\sigma=0$, $\rho=-\overline{\rho}$ and, by a choice of frame, $\epsilon=0$.  By the identity (\ref{e:id1}) $D\rho=0$ and projecting the curvature condition (\ref{e:sas}) onto the orthonormal frame we find that it reduces to the single condition $\rho^2=-1$.

Thus  ($M$, $g$) is Sasakian if and only if it admits a divergencefree, shearfree geodesic of constant twist.

$\Box$

In such a case, writing $\rho= i$, we have that
\bee\label{e:sasric1}
R_{00}=2 \quad R_{0-}= D\tau -\tau i \quad R_{+-}= \delta\tau+\overline{\delta}\;\overline{\tau}  -2\tau\overline{\tau} \quad R_{++}= R_{0+}=0
\ee
\bee\label{e:sasric2}
{\shf}R= \delta\tau+\overline{\delta}\overline{\tau} -2\tau\overline{\tau}+1  .
\ee
and the Bianchi identity (\ref{e:bian1}) reduces to 
\[
D(\delta\tau+\overline{\delta}\overline{\tau}-2\tau\overline{\tau})=0.
\]
Thus we conclude that $DR=0$, that is, the scalar curvature is constant along the geodesic congruences.

The equations then to be solved are

\bee \label{e:ell1}
D\tau = \tau i
\ee
\bee \label{e:ell2}
\delta\tau+\overline{\delta}\;\overline{\tau}=2\tau\overline{\tau}-1+{\hf}R
\ee
\bee \label{e:comm1ella}
D\Omega=-i \Omega
\ee
\bee \label{e:comm1ellb}
D\xi^a=-i \xi^a
\ee
\bee \label{e:comm2ella}
\delta\overline{\Omega}-\overline{\delta}\Omega=\overline{\tau}\;\overline{\Omega}-\tau\Omega -2  i
\ee
\bee \label{e:comm2ellb}
\delta\overline{\xi}^a-\overline{\delta}\xi^a=\overline{\tau}\;\overline{\xi}^a-\tau\xi^a.
\ee
The solutions to equations (\ref{e:ell1}), (\ref{e:comm1ella}) and (\ref{e:comm1ellb}) are
\[
\tau=\tau_{0}\;e^{i r}
\]
\[
\Omega=\Omega_{0}\;e^{-i r}
\]
\[
\xi^a=\xi_{0}^a\;e^{-i r}.
\]
where a subscript 0 indicates independance from $r$.

We can choose coordinates ($x^2,x^3$) so that 
\[
\xi^a_{0}=P_{0}(\delta^a_2+i\delta^a_3),
\]
and a frame so that $P_{0}$ is a real non-zero function of ($x^2,x^3$).  The remaining equations (\ref{e:ell2}), (\ref{e:comm2ella}) and (\ref{e:comm2ellb}) become
\bee\label{e:finell1}
2P_{0}\left[ \frac{\partial \tau_{0}}{\partial\overline{z}}+  \frac{\partial \overline{\tau}_{0}}{\partial z}\right] +(\Omega_{0}\tau_{0}-\overline{\Omega}_{0}\;\overline{\tau}_{0}) i=2\tau_{0}\overline{\tau}_{0}-1+{\hf}R
\ee
\bee\label{e:finell2}
2P_{0}\left[ \frac{\partial \overline{\Omega}_{0}}{\partial\overline{z}} -  \frac{\partial \Omega_{0}}{\partial z}\right] + 2\Omega_{0}\;\overline{\Omega}_{0}  i =\overline{\tau}_{0}\;\overline{\Omega}_{0}-\tau_{0}\Omega_{0}-2i
\ee
\bee\label{e:finell3}
2\frac{\partial P_{0}}{\partial \overline{z}}=\overline{\tau}_{0}-\Omega_{0} i,
\ee
where we have introduced the complex coordinates $z=x^2+ix^3$ and $\overline{z}=x^2-ix^3$ and
\[
\frac{\partial}{\partial z}=\frac{1}{2}\left(\frac{\partial}{\partial x^2}-i\frac{\partial}{\partial x^3}\right)
\qquad\qquad
\frac{\partial}{\partial \overline{z}}=\frac{1}{2}\left(\frac{\partial}{\partial x^2}+ i\frac{\partial}{\partial x^3}\right).
\]
  Now, differentiating (\ref{e:finell3}) and using (\ref{e:finell1}) and (\ref{e:finell2}) yields
\bee\label{e:tw}
8P_{0}^2\frac{\partial^2 \ln P_{0}}{\partial z \partial \overline{z}}=1+{\hf}R.
\ee

A relabelling of the rays $r\rightarrow r'=r+f(x^2,x^3)$ allows us to set the imaginary part of $\Omega_{0}$ to zero.  The final equation to be solved (equation (\ref{e:finell2})) is
\[
P_{0}^2\left[\frac{\partial}{\partial \overline{z}}\left(\frac{\Omega_{0}}{P_{0}}\right) - \frac{\partial}{\partial z}\left(\frac{\Omega_{0}}{P_{0}}\right)\right]=-  i
\]
or, if ($x^2=u$, $x^3=v$)
\bee\label{e:last2}
P_{0}^2\frac{\partial}{\partial v}\left(\frac{\Omega_{0}}{P_{0}}\right)=-  1.
\ee
Thus we have the following local moduli theorem for Sasakian structures:

\vspace{0.2in}

\noindent{\bf Theorem 3}

{\it
For any Sasakian structure on $R^3$ there exists local coordinates ($r,u,v$) such that the Killing vector is $\frac{\partial}{\partial r}$ and the metric is
\[
ds^2=\left[dr+\int\frac{1}{P_0^2}dv\;.\;du\right]^2 +\frac{1}{2P_0^2}\left[du^2+dv^2\right],
\]
where $P_0$ is an arbitrary nowhere zero function of $u$ and $v$. The arbitrary function of integration can be removed by a relabelling of $r$.
}

The question of uniqueness is dealt with by:

\vspace{0.2in}

\noindent{\bf Theorem 4}

{\it
Two Sasakian structures  given by
\[
ds^2=\left[dr+\int\frac{1}{P_0^2}dv\;.\;du\right]^2 +\frac{1}{2P_0^2}\left[du^2+dv^2\right]
\]
\[
d\tilde{s}^2=\left[d\tilde{r}+\int\frac{1}{\tilde{P}_0^2}d\tilde{v}\;.\;d\tilde{u}\right]^2 +\frac{1}{2\tilde{P}_0^2}\left[d\tilde{u}^2+d\tilde{v}^2\right],
\]
are contact isometric if and only if 
\bee\label{e:anal}
\tilde{P}_0^2\left(\frac{d z}{d w}\right)\overline{\left(\frac{d z}{d w}\right)}=P_0^2,
\ee
where $z=u+iv$ is an analytic function of $w=\tilde{u}+i\tilde{v}$.
}

\noindent{\bf Proof}:

Here, two smooth riemannian manifolds ($M_1$, $g_1$) and ($M_2$, $g_2$) with contact structures given by contact 1-forms $\alpha_1$ and $\alpha_2$, respectively, are {\it contact isometric} if there exists a diffeomorphism $\phi:M_1\rightarrow M_2$ such that $\phi^*g_2=g_1$ and $\phi^*\alpha_2=\alpha_1$.  

Suppose then that the two Sasakian structures are contact isometric.  Then 
\[
dr+\int\frac{1}{P_0^2}dv\;.\;du = d\tilde{r}+\int\frac{1}{\tilde{P}_0^2}d\tilde{v}\;.\;d\tilde{u} ,
\]
and 
\[
\frac{1}{P_0^2}\left[du^2+dv^2\right] =
\frac{1}{\tilde{P}_0^2}\left[d\tilde{u}^2+d\tilde{v}^2\right].
\]
From these we see that
\[
r=\tilde{r}+f(\tilde{u},\tilde{v})\qquad\qquad u=u(\tilde{u},\tilde{v}) \qquad\qquad v=v(\tilde{u},\tilde{v}).
\]
Moreover, we must have $z=z(w)$ where $z=u+iv$ is an analytic function of $w=\tilde{u}+i\tilde{v}$, and $P_0$ and $\tilde{P}_0$ are related by (\ref{e:anal}), as required.

Conversely, suppose (\ref{e:anal}) holds.  Then by Theorem 3 there must be a relabelling of the rays $ r=\tilde{r}+f(w,\overline{w})$ such that the transformation
\[
r=\tilde{r}+f(w,\overline{w})\qquad\qquad z=z(w) \qquad\qquad \overline{z}=\overline{z(w)}
\]
takes the first Sasakian structure to the second.

$\Box$

\vspace{0.2in}

In addition, we have the following

\vspace{0.2in}

\noindent {\bf Corollary 5}

{\it
Given any smooth function $R$ of $u$ and $v$, there exists a Sasakian structure on $R^3$ such that $R$ is the scalar curvature of the associated metric.
}

\noindent{\bf Proof:}

By a theorem of Kazdan and Warner \cite{kaw}, equation (\ref{e:tw}) has a non-zero solution $P_0$ for a given smooth function $R(u,v)$. Then, integrating equation (\ref{e:last2}) gives the required Sasakian structure.

$\Box$

\section{ Sasakian $\eta$-Einstein Structures}

A Sasakian structure is said to be $\eta$-Einstein if there exist constants $a$ and $b$ such that
\bee\label{e:saseta}
R_{ij}=a\;g_{ij}+b\;\eta_i\eta_j,
\ee
where $\eta$ is the 1-form associated with the unit Killing vector.  For such a structure on a 3-manifold the scalar curvature is constant 
\[
R=3a+b.
\]
Indeed, it is not hard to see that a constant scalar curvature Sasakian 3-manifold must be $\eta$-Einstein.  Projecting (\ref{e:saseta}) onto the Killing vector direction, using equation (\ref{e:sasric1}), we find that $a+b=2$, and so $R=2(a+1)$.

The Tanaka-Webster curvature of a Sasakian structure on a 3-manifold is given by \cite{guil}
\[
W=\frac{1}{4}(R+2),
\]
and so a constant scalar curvature Sasakian 3-manifold has constant Tanaka-Webster curvature.  For a Sasakian $\eta$-Einstein manifold this simplifies to
\[
W=\frac{1}{2}(a+2).
\]
In order to find all such structures we start with equation (\ref{e:tw}):
\bee\label{e:tw2}
4P_{0}^2\frac{\partial^2 \ln P_{0}}{\partial z \partial \overline{z}}=W.
\ee
To solve this, we need to consider the cases $W>0$, $W=0$ and $W<0$ separately.

\subsection{Sasakian structures with constant $W>0$}

In an appropriate coordinate system, the solution to (\ref{e:tw2}) is
\[
P_{0}={\hf}\sqrt{W}(1+z\overline{z}),
\]
while equation (\ref{e:last2}) yields
\[
\Omega_{0}=-\frac{1}{\sqrt{W}} \left[\frac{v}{1+u^2}+\frac{(1+u^2+v^2)}{(1+u^2)^{3/2}}
    \tan^{-1}\left(\frac{v}{\sqrt{1+u^2}}\right)\right].
\]
In fact, the metric can be much simplified by the transformation ($r,u,v$)$\rightarrow$($\rho,\theta,\phi$) to Euler coordinates given by
\[
r =\frac{\rho +\phi}{W}-\frac{2\cos\phi\tan({\sdv{\theta}{2}})}{W\sqrt{1+\cos^2\phi\tan^2({\sdv{\theta}{2}}})}\tan^{-1}\left[\frac{\sin\phi\tan({\sdv{\theta}{2}})}{\sqrt{1+\cos^2\phi\tan^2({\sdv{\theta}{2}})}}\right]
\]
\[
u=\cos\phi\tan({\sdv{\theta}{2}})
\]
\[
v=\sin\phi\tan({\sdv{\theta}{2}}).
\]
The metric reduces to the form:
\[
ds^2=\frac{1}{2W}\left[d\theta^2+\sin^2\theta d\phi^2+\frac{2}{W}(d\rho+\cos\theta d\phi)^2\right],
\]
and the contact 1-form becomes:
\[
\alpha=\frac{1}{W}(d\rho+\cos\theta d\phi).
\]
The special case $W=2$ gives the standard Sasakian structure associated with the round metric on $S^3$ (see \cite{guil}), while for $W\neq2$ the metric is not Einstein.  Remarkably, for $0<W<2$ the metric is homothetic to the Berger sphere \cite{tod}.  This metric can be obtained from the round metric on $S^3$ by deforming the metric along the fibres of the Hopf fibering, which form the shearfree geodesic congruences.  This is an Einstein-Weyl space and the Weyl 1-form is just a constant multiple of the contact 1-form $\alpha$.

\subsection{Sasakian structures with $W=0$}

In this case, the solution to (\ref{e:tw2}) can be set to
\[
P_{0}=\frac{1}{\sqrt{2}},
\]
while equation (\ref{e:last2}) yields
\[
\Omega_{0}=-\sqrt{2}v.
\]
Thus the metric is
\[
ds^2=\left[dr+2vdu\right]^2+du^2+dv^2.
\]
which has negative scalar curvature $R=-2$.  This is homothetic to the standard Sasakian structure on the 3-manifold $Nil$.

\subsection{Sasakian structures with constant $W<0$}

An appropriate choice of coordinate gives the solution to (\ref{e:tw2}) as
\[
P_{0}={\hf}\sqrt{-W}(1-z\overline{z}),
\]
while equation (\ref{e:last2}) now yields
\[
\Omega_{0}=-\frac{1}{\sqrt{-W}} \left[\frac{v}{1-u^2}+\frac{(1-u^2-v^2)}{(1-u^2)^{3/2}}
    \tanh^{-1}\left(\frac{v}{\sqrt{1-u^2}}\right)\right].
\]
The metric can be much simplified by the transformation ($r,u,v$)$\rightarrow$($\rho,\theta,\phi$) given by
\[
r =\frac{\rho +\phi}{W}+\frac{2\cos\phi\tanh({\sdv{\theta}{2}})}{W\sqrt{1-\cos^2\phi\tanh^2({\sdv{\theta}{2}}})}\tanh^{-1}\left[\frac{\sin\phi\tanh({\sdv{\theta}{2}})}{\sqrt{1-\cos^2\phi\tanh^2({\sdv{\theta}{2}})}}\right]
\]
\[
u=\cos\phi\tanh({\sdv{\theta}{2}})
\]
\[
v=\sin\phi\tanh({\sdv{\theta}{2}}).
\]
The metric reduces to the form:
\[
ds^2=-\frac{1}{2W}\left[d\theta^2+\sinh^2\theta d\phi^2-\frac{2}{W}(d\rho+\cosh\theta d\phi)^2\right],
\]
and the contact 1-form becomes:
\[
\alpha=\frac{1}{W}(d\rho+\cosh\theta d\phi).
\]
For $W=-2$ this is the locally homogenous metric on $\tilde{SL_2R}$

\section{Conformally Flat Sasakian Structures}

In this section we prove the following theorem:

\vspace{0.2in}

\noindent{\bf Theorem 6}

{\it
A conformally flat Sasakian 3-manifold is covered by $S^3$ with the standard structure.
}

\noindent{\bf Proof}:

A metric is conformally flat in three dimensions iff the Weyl-Schouten tensor vanishes:
\bee\label{e:conff}
\nabla_{[i}C_{j]k}=0,
\ee
where
\[
C_{ij}=R_{ij}-{\dv{1}{4}}Rg_{ij}.
\]
For a Sasakian structure the non-zero components of $C_{ij}$ are
\[
C_{00}=-{\hf}(\delta\tau+\overline{\delta}\overline{\tau}-2\tau\overline{\tau}-3)
\]
\[
C_{+-}={\hf}(\delta\tau+\overline{\delta}\overline{\tau}-2\tau\overline{\tau}-1).
\]
Projecting equation (\ref{e:conff}) onto the frame  adapted to the Sasakian structure we find that
\[
\delta C_{00}=0
\]
\[
C_{00}-C_{+-}=0
\]
\[
DC_{+-}+2i(C_{00}-C_{+-})=0
\]
\[
\delta C_{+-}=0.
\]
From these we conclude that $C_{00}=C_{+-}={\hf}$, $R=6$ and finally
\[
R_{ij}={\dv{1}{3}}Rg_{ij}.
\]
Thus, the metric is Einstein.  Geiges \cite{geiges} has shown that if the manifold is closed then it is diffeomorphic to $S^3/\Gamma$.  If it is open it has been shown \cite{guil} that it is locally covered by $S^3$.

$\Box$

\section{Remarks}
Consider a closed 3-manifold $M$ with a regular Sasakian structure, that is, the integral curves of the Killing vector $e_0$ are generated by a free $S^1$ action on $M$.  Then $M$ fibres over a riemann surface $\Sigma$.  Since $e_0$ is Killing, the metric  on $\Sigma$ is just the induced metric from $M$, while the contact form $\theta^0$ is a connection in the principal bunble $M\rightarrow \Sigma$.   The curvature of the $\Sigma$ is the Tanaka-Webster curvature $W$, so the Euler characteristic of $\Sigma$ is 
\[
\chi=\frac{1}{2\pi}\int_\Sigma W dA
\]
On the other hand the Euler number of the bundle can be defined in terms of the curvature of the connection $\theta^0$
\[ 
e=\int_\Sigma F=\int_\Sigma d\theta^0=\lambda\int_\Sigma \theta^1\wedge\theta^2\neq 0
\]
Amongst the eight locally homogenous geometries, the ones that can admit  Sasakian structures are seifert fibred spaces with $e\neq 0$, which turn out to be

\begin{center}
\vskip .2truein
\begin{tabular}{|c||c|c|c|} \hline
Euler Characteristic&$\chi>0$&$\chi=0$ &  $\chi<0$ \\ \hline
Manifold& $S^3$  &$Nil$ &$\tilde{SL_2R}$    \\ \hline
\end{tabular}
\end{center} 

Theorems 3 and 4 essentially say that, locally, all Sasakian structures come from a line bundle over a riemann surface (a Kaluza-Klein construction) and that the moduli space of such structures corresponds to the riemann structure on the surface.  When $W$ is constant we have found the complete set of Sasakian structures which includes the three locally homogenous geometries as listed above.  In addition, we have found many examples of Sasakian structures with negative scalar curvature:

\begin{center}
{\bf Sasakian $\eta$-Einstein 3-manifolds }
\vskip .2truein
\begin{tabular}{|c||c|c|c|} \hline
  & $R>0$ & $R=0$ & $R<0$  \\ \hline\hline
$W>0$&$W>{\shf}$ &$W={\shf}$ & $0<W<{\shf}$ \\ \hline
$W=0$&None &None &  All \\ \hline
$W<0$&None &None &  All \\ \hline
\end{tabular}
\end{center}

Finally, the function $P_0$ is related to the recently discovered \cite{gkn} Sasakian potential $K$ by
\[
\frac{1}{P_0^2}=4\frac{\partial K}{\partial z \partial \overline{z}},
\]
although we have used different coordinate conditions for the description of the contact 1-form.


\begin{references}
\bibitem{sas} S. Sasaki, T\^{o}hoku Math. J. {\bf 2} 459 (1960).
\bibitem{bag1} C. P. Boyer and K. Galicki {\it 3-Sasakian manifolds} Preprint  MPI-98/19 (1998).
\bibitem{bag2} C. P. Boyer and K. Galicki {\it On Sasakian-Einstein geometry} UNM Preprint  (1998).
\bibitem{geiges} H. Geiges, T\^{o}hoku Math. J. {\bf 49} 415 (1997).
\bibitem{belgun} F. Belgun {\it Normal CR structures on compact 3-manifolds } Preprint  (2000).
\bibitem{guil} B. S. Guilfoyle,{\it Einstein metrics adapted to contact structures on 3-manifolds} Preprint  http://arXiv.org/abs/math/0012027 (2000).
\bibitem{tanno2} S. Tanno, Proc. Japan Acad. {\bf 43 } 581 (1967).
\bibitem{kaw} J. L. Kazdan and F. W. Warner, Ann. of Math. {\bf 99} 203 (1974).
\bibitem{tod} K. P. Tod, J. London Math. Soc. {\bf 45 } 341 (1992).
\bibitem{gkn} M. Godlinski, W. Kopczynski and P. Nurowski, {\it Locally Sasakian manifolds } Preprint  http://arXiv.org/abs/math/0005072 (2000).
\end{references}
\end{document}